# Equation-Free Multiscale Computational Analysis of Individual-Based Epidemic Dynamics on Networks

Constantinos I. Siettos

School of Applied Mathematics and Physical Sciences,

National Technical University of Athens, GR 157 80, Athens, Greece

**Abstract**

The surveillance, analysis and ultimately the efficient long-term prediction and control of epidemic dynamics appear to be one of the major challenges nowadays. Detailed atomistic mathematical models play an important role towards this aim. In this work it is shown how one can exploit the Equation Free approach and optimization methods such as Simulated Annealing to bridge detailed individual-based epidemic simulation with coarse-grained, systems-level, analysis. The methodology provides a systematic approach for analyzing the parametric behavior of complex/ multi-scale epidemic simulators much more efficiently than simply simulating forward in time. It is shown how steady state and (if required) time-dependent computations, stability computations, as well as continuation and numerical bifurcation analysis can be performed in a straightforward manner. The approach is illustrated through a simple individual-based epidemic model deploying on a random regular connected graph. Using the individual-based microscopic simulator as a black box coarse-grained timestepper and with the aid of Simulated Annealing I compute the coarse-grained equilibrium bifurcation diagram and analyze the stability of the stationary states sidestepping the necessity of obtaining explicit closures at the macroscopic level under a pairwise representation perspective.

## 1. Introduction

No doubt, the history of mankind has been shaped by the pitiless ravages of pandemics. Whole nations and civilizations have been wiped off the map through the ages. The list is long: biblical pharaonic plagues which hit Ancient Egypt in the middle of Bronze Age around 1715 B.C. [1], the "λοιμός" –the unidentified plague- which stroke Athens from 430 to 425 B.C. and set the end of the Periclean golden era [2], the "cocoliztli" epidemics occurred during the 16th century resulting to some 13 million deaths decimating the Mesoamerican native population [3], the Black Death Bubonic Plague which burst in Europe in 1348 and is estimated to have killed over 25 million people in just five years [4]. Ninety years ago the pandemic influenza virus of 1918–1919 swept through America, Europe, Asia and Africa smashing the globe: the death-toll was around 40 million people. Two one-year less severe influenza pandemics followed in the next decades: the 1957 and the 1963 influenza pandemics resulted to two and one million deaths respectively [5].

Forty four years later even though there has been an immense progress in combating infections, history and studies warn us that complex multiscale interactions between a host of factors ranging from the micro host-pathogen and individual-scale host-host interactions to macro-scale ecological, social, economical and demographical conditions across the globe aggravated by the well-fare and heath-care system degradation in the poorest parts of the world result to the re-emergence of latent as well as the appearance of newly emergent diseases.

In June 2009 swine flu was declared as a level six pandemic by the World Health Organization (WHO). The Center for Disease Control and Prevention of the USA estimated that between 41 million and 84 million cases of 2009 H1N1 flu occurred in the USA between April 2009 and January 16, 2010 [6]. WHO reported that as of 24 January 2010 the swine-flu pandemic (lab-confirmed) has caused the death of more than 15,300 world-wide across 212 countries [7]. Even though the consequences have been relatively mild compared to the seasonal flu it has already been one for the ages. The most worrisome: the possibility of the emergence of a more dangerous pandemic wave in the near future cannot be excluded.

The critical question(s) is(are) not whether a new pandemic will arise but when it will, how it is going to spread, how deadly it will be, who should get the vaccine when not all can, how likely are multiple waves of re-emergence and what type of intervention may be applied to stop the spread. Unfortunately, even with all the advances, we still don't have robust answers. The breaking news of the WHO reminds us of our vulnerability.

Mathematical models and systems theory are playing a most valuable role in shedding light on the problem and for helping make decisions. The studies have proceeded mainly on two fronts. On the one hand, they are the "continuum models" describing the coarse-grained dynamics of the disease in the population (see e.g. [8-11]) one might, for example study a model for the evolution of the disease as a distributed function of the age and the time since vaccination (see e.g. [12]). Such models can be explored using powerful analysis techniques for ordinary and partial differential equations. However, due to the complexity of the phenomena, available continuum models are often qualitative caricatures of the reality. On the other hand there are the so-called object-based models where the complex system is viewed as a network of interacting discrete entities (individuals or objects). This group includes models ranging from cellular automata (see e.g. [13-16]) to individual-based [17-19] and very detailed agent-based models on various social network graphs [20-23].

But while object-based models are becoming the tools of choice of today, they are only half the battle. At the end of the day there is a basic question to be answered: how much could we trust the outcomes in a real crisis? Such state-of-the-art models are -as all models- just approximations of the real system they aspire to represent. Due to the inherent extraordinarily complexity of the problem, they are built with incomplete knowledge and for that reason they are flashing a "note of caution" on parameter and rule inaccuracies.

To date what is usually done with the detailed object oriented descriptions lacking explicit macroscopic descriptions is simple simulation: set up many initial conditions, for each initial condition create a large enough number of ensemble realizations, probably change some of the rules and then run the detailed dynamics for a long time to investigate how things like vaccination policies, malignancy of a virus -as this may be expressed in terms of the reproduction number-, and resource availability may influence the spread of an outbreak. However, this simple simulation is most of the times inadequate for the systematic analysis, control policies design and optimization of epidemic dynamics.

One of the most critical, yet unresolved issues in the area is this: understand how to systematically bridge the gap between the micro-scale, where detailed bio-information on the immune mechanisms, host-virus and host-host interactions is often available, and the city or country scale, where the disease emerges, the strategic combat-policy questions arise and the answers are required. However due to the emerging complex, nonlinear, stochastic nature of the dynamics of such models, and the intrinsic multiplicity of scales at which the relevant objects interact, makes the systematic analysis and design at the macroscopic level an overwhelming task: good macroscopic, evolution equations in closed form, that will allow us to predict- in a systematic way- the evolution under various scenarios and design control/ depression policies, cannot be written in a straightforward manner.

To deal with the above problem, various moment-approximations such as mean field and pair-wise formulations [18,24-26] have been proposed in order to extract analytical closed macroscopic models of the underlying detailed dynamics. These, try to relate higher-order moments of the atomistic interactions- causing the spread of the disease-evolving over a social network to a few, low-order ones. For example, pair-approximation schemes try to involve higher–order moments by analytically writing down moment closures which usually relate densities of triples of states of individuals (corresponding the third-order moments) with densities of pairs of individuals (corresponding to the second order moments) over the network [27-29]. However, while these approaches offer a good starting point for analyzing network-based epidemic dynamics in a more formal way, they introduce systematic bias and therefore may miss important qualitative characteristics at the coarse-grained level as the interaction dynamics become more complex (see for example the discussion in [18, 30]).

New computational methodologies that could systematically extract coarse-grained, emergent dynamical information by bridging complex individual-based modeling with macroscopic, systems-level, continuum numerical analysis, control and optimization methods resolving in a systematic way the above problems, would facilitate the exploration and therefore have the potential significantly contribute to our understanding, predicting and designing – of better public health strategies to combat emergent epidemics.

Over the past few years it has been demonstrated that the so called "Equation-Free" approach can be used to establish the missing link between traditional numerical analysis, control design tools and microscopic/ stochastic simulators [31-39]. This mathematical-assisted approach serves an on-demand identification-based computational approach enabling microscopic/ atomistic-based models to perform system-level tasks bypassing the necessity of deriving good *explicit* closures.

Here it is shown how one can exploit the approach with the help of Simulated Annealing (SA) [40-42] to study in a computational strict, systematic and efficient way the macroscopic, emergent behavior of individual-based epidemic models on networks within the pairwise representation context. For illustration purposes, I analyse the coarse-grained dynamics of a simple individual-based epidemic model deploying on a fixed random regular network (RRN) serving as a caricature of the underlying social/interaction structure. The model describes the spread of a hypothetical infectious disease in a population of constant size. Here each individual can be in one of three states: susceptible, infected or recovered and interacts with four other individuals. The states of the individuals change over discrete time in a probabilistic manner according to simple rules involving their own states and the states of their links. The relatively simple rules governing the interactions between the individuals result to a fundamental feature of such problems which is the emergence of complex dynamics in the coarse-grained level such as the multiplicity of coarse-grained stationary states leading to hysteresis phenomena (see e.g. [43]); more complex behaviour such as recurrent situations and chaos have also been observed in real-world epidemics [44-46]. Using the atomistic simulator as a black-box, I construct the coarse-grained bifurcation diagram and perform stability analysis of the computed solutions.

To this end, the paper is organized as follows: in the next section I describe the kernel of the Equation-free approach and present how this can be combined with SA for multiscale calculations on epidemic networks. In section three I describe the individual-based epidemic model while in section four I illustrate the results of the analysis. Finally the main conclusions of this work are summarized in section five.

## 2. The Computational framework for multiscale computations on epidemic simulators

### 2.1 A short description of The Equation-Free concept

For detailed individual-oriented simulators of epidemics, the closures required to extract good representative models in the continuum/macroscopic level are -due to the inherent multiscale complexity and heterogeneity - most of the times not available and rather difficult to derive. Besides, simulation forward in time is but the first of systems level computational tasks one wants to explore while analyzing the parametric behavior of such large-scale individualistic epidemic models. When this is so, the so-called Equation-free approach to modeling and multiscale system analysis, a novel computational framework, can be used to circumvent the need for obtaining an explicit continuum model in closed form [31-39, 47-51]. Indeed, through appropriate computations of the detailed models, one is able to estimate the same information that a continuum model would allow us to compute from an explicit formula. Using this framework, steady state and stability computations, as well as continuation and numerical bifurcation analysis of the complex-emergent dynamics can be performed in a full-computational manner bypassing the need of analytical derivation of closures for the macroscopic-level equations.

The main assumption behind the methodology is that a coarse-grained model for the dynamics at the macroscopic/continuum level exists and closes in terms of a few coarse-grained variables which are usually the low-order moments of the microscopically evolving distributions and simultaneously the apparent observables of the evolving phenomenon (e.g. the rate of spread, the distribution of the diseases within the population as a function of the age).

Let us start the presentation of the methodology supposing that we do not have available the explicit macroscopic equations in a closed form, but we do have an evolving microscopic (very) large scale computational model.

Given the (microscopic) distribution of the system $U_k \equiv U(t_k) \in R^N$, $N >> 1$ at time $t_k = kT$ the detailed simulator reports the values of the state variables after a time interval $T$, i.e.:

$$U_{k+1} = \wp_T(U_k, p), \quad (1)$$

where $\wp_T : R^N \times R^m \to R^N$ is the time-evolution operator, $p \in R^m$ is the vector of systems parameters.

The main assumption behind the coarse timestepper concept, the core of the Equation-Free approach, is that a coarse-grained model for the fine-scale dynamics (1) exists and closes in terms of a few coarse-grained variables, say, $x \in R^n, n << N$.

Usually, these are the low-order moments of the microscopically evolving distributions. The existence of a coarse-grained model implies that the higher order moments, say, $y \in R^{N-n}$, of the distribution $U$ become, relatively fast over the coarse time scales of interest, functions of the few lower ones, $x$. This scale separation could be viewed in the form of a singularly perturbed system of the following form:

$$x_{k+1} = h_s(x_k, y_k, \varepsilon, p) \tag{2a}$$

$$\varepsilon\, y_{k+1} = h_f(x_k, y_k, \varepsilon, p) \tag{2b}$$

with $\varepsilon$ a sufficiently small real number. Eq. (2a) corresponds to the "slow" coarse-grained dynamics while Eq. (2b) to the "fast" ones.

Let us suppose that the following assumptions also hold:

(i) the functions $h_s$, $h_f$ are of class $C^r, r \geq 2$,

(ii) we can apply the implicit function theorem to find a function relating the coarse-grained fast and slow dynamics when $\varepsilon = 0$ in the form of

$$y = q(x, p, \varepsilon = 0) \tag{3a}$$

such that

$$0 = h_f(x_k, q(\cdot), 0, p) \tag{3b}$$

The above assumption guarantees the existence of a slow manifold $M_0$ defined by (3a). On $M_0$, the dynamics of the original system given by (2) can be described by the reduced system

$$x_{k+1} = h_s(x_k, q(\cdot), 0, p) \tag{4}$$

(iii) The manifold $M_0$ is hyperbolic, i.e. the Jacobian matrix $\nabla h_f(\cdot)$ is not singular. Hence, all the eigenvalues of $\nabla h_f(\cdot)$ evaluated $\forall (x, y) \in M_0$ have either positive or negative real parts. This implies that the fast variables $y$ can be in principle decomposed into stable and unstable subspaces corresponding to attracting and repelling manifolds respectively. From now on we will assume that

$$Re\, \lambda_i \{\nabla h_f(\cdot)\} < 0, i = 1, 2, ...N\text{-}n \tag{5}$$

The above basic assumption ensures that the fast dynamics of (2) will converge to their quasi-steady state (3a) and will not drive the system towards infinity [52, 53]. Hence, a well defined coarse-grained low-order model in principle exists and its dynamics can be approximated by the reduced model (4).

Fenichel's theorem [54] can be used at this point to reduce the overall dynamics defined by Eq. (2) to the following system:

$$x_{k+1} = h_s(x, q(x, p, \varepsilon), p) \tag{6a}$$

where

$$y = q(x, p) \tag{6b}$$

determines the slow manifold defined by

$$M_e = \{(\mathbf{x}, \mathbf{y}) \in R^n \times R^{N-n} : \mathbf{y} = q(\mathbf{x}, \mathbf{p})\} \tag{7}$$

on which the coarse-grained dynamics of the system evolve after a fast transient phase (see figure 1). $M_\varepsilon$ is $O(\varepsilon)$ close to the equilibrium manifold $M_0$, i.e. $M_\varepsilon \to M_0$, $\varepsilon \to 0$.

What the methodology does, in fact, is providing a closure such as (6b) "on demand" in a computational manner, without writing it down.

In a nutshell, the computation methodology consists of the following steps (see also figure 2):

**(a)** Choose the coarse-grained statistics of interest for describing the long-term behavior of the system and an appropriate representation for them (for example the densities of the susceptible and infected individuals in the population).
**(b)** Choose an appropriate lifting operator $\mu$ from the continuum description $\mathbf{x}$ to the individual-based description $\mathbf{U}$ on the network. For example, $\mu$ could make random susceptible and infection assignments over the network consistent with the respective densities.
**(c)** Prescribe a continuum initial condition at a time $t_k$: $\mathbf{x}_{t_k}$.
**(d)** Transform this initial condition through lifting to one (or more) consistent individual-based realization(s) $\mathbf{U}_{t_k} = \mu \mathbf{x}_{t_k}$.
**(e)** Evolve thi(e)s(e) realization(s) using the microscopic (individual-based) model for a desired time $T$, generating the $\mathbf{U}_{t_{k+1}}$, where $t_k = kT$.
**(f)** Obtain the restrictions $\mathbf{x}_{t_{k+1}} = \aleph \mathbf{U}_{t_k}$.

This, constitutes the *coarse timestepper*, or *coarse time-T map* that, given an initial coare-grained state of the system $\mathbf{x}_{t_k}$, $\mathbf{p}$ at time $t_k$ will report the result of the integration of the individual-based rules after a given time-horizon $T$ (at time $t_{k+1}$), i.e.

$$\mathbf{x}_{t_{k+1}} = \mathbf{\Phi}_T(\mathbf{x}_{t_k}, \mathbf{p}), \tag{8}$$

where $\mathbf{\Phi}_T: R^n \times R^m \to R^n$ having $\mathbf{x}_k$ as initial condition.

The assumption of the existence of a coarse-grained the temporal evolution operator $\mathbf{\Phi}_{T_h}$, which is assumed to be unavailable analytically in closed form, implies that the higher order moments of the distributions become, relatively quickly over the coarse-time scales of interest, "slaved" to the lower, few, "master" ones.

At this point one can wrap around the coarse-grained input-output map (8) a fixed point iterative scheme order to compute fixed point or periodic solutions at certain values of the parameter space.
For example for low-order systems coarse-grained equilibria can be obtained as fixed points, of the map $\mathbf{\Phi}_T$:

$$\mathbf{x} - \mathbf{\Phi}_T(\mathbf{x}, \mathbf{p}) = \mathbf{0} \tag{9}$$

using a simple Newton-Raphson scheme.
The tracing of the branches in a one-dimensional parameter space through regular turning points can be achieved using a continuation technique such as the linearized pseudo arc-length procedure, augmenting the systems equation with the condition [55]:

$$N(\mathbf{x}, p) = \mathbf{\alpha} \cdot (\mathbf{x} - \mathbf{x}_1) + \beta(p - p_1) - \delta s = 0, \tag{10a}$$

where

$$\boldsymbol{\alpha} \equiv \frac{(\boldsymbol{x}_1 - \boldsymbol{x}_0)^T}{\delta s}, \beta \equiv \frac{(p_1 - p_0)}{\delta s}, \tag{10b}$$

and $\delta s$ is the pseudo arc-length continuation step. Here $(\boldsymbol{x}_0, p_0)$ and $(\boldsymbol{x}_1, p_1)$ are two already computed solutions. Equation (6) constrains the "next" equilibrium of (9) to lie on a hyperplane perpendicular to the tangent of the bifurcation diagram at $(\boldsymbol{x}_1, p_1)$, approximated through $(\boldsymbol{a}, \beta)$ and at a distance δs from it. Now, the computation of the "next" fixed point involves the iterative solution of the following linearized system:

$$\begin{bmatrix} \boldsymbol{I} - \dfrac{\partial \boldsymbol{\Phi}_T}{\partial \boldsymbol{x}} & -\dfrac{\partial \boldsymbol{\Phi}_T}{\partial p} \\ \boldsymbol{a} & \beta \end{bmatrix} \begin{bmatrix} \delta \boldsymbol{x} \\ \delta p \end{bmatrix} = -\begin{bmatrix} \boldsymbol{x} - \boldsymbol{\Phi}_T(\boldsymbol{x}, p) \\ N(\boldsymbol{x}, p) \end{bmatrix} \tag{11}$$

The local stability of the above system around equilibria is determined by the the eigenvalues of the Jacobian matrix $\dfrac{\partial \boldsymbol{\Phi}_T}{\partial \boldsymbol{x}}$.

As the size *n* of the coarse-grained system increases one can use matrix-free iterative methods such as Newton-GMRES [56] and Arnoldi's procedure [57] in order to find fixed points of (9) and approximate the most critical eigenmodes that determine the stability of the solutions, respectively.

If a periodic oscillatory behaviour with a period of $\Delta T$ is observed then one seeks for solutions which satisfy

$$\boldsymbol{x}(t) = \boldsymbol{x}(t + \Delta T), \tag{12}$$

In this case periodic solutions can be computed as fixed points of the map

$$\boldsymbol{x} - \boldsymbol{\Phi}_{\Delta T}(\boldsymbol{x}, p) = \boldsymbol{0} \tag{13}$$

The unknown period can be computed by augmenting the above system of equations by the so-called phase constraint (also called a pinning condition)

$$g(\boldsymbol{x}, p, \Delta T) = 0, \tag{14}$$

which factors out the infinite members of the family of periodic solutions in (12) [58].

**2.2 Driving the epidemic simulator to the slow manifold**

**2.2.1 The slow and fast variables of epidemic models within networks**

For the scheme to be accurate, the overall procedure has to be applied when the system evolves on the slow manifold (i.e., the fixed point iteration (3) has to be solved *on* the slow manifold). If the time required for trajectories emanating from initial conditions off the slow manifold to reach the slow manifold is very small compared to *T* (i.e. when the microscopic initial conditions are close enough to the slow manifold), the above requirement is satisfied for any practical means. In general, we would expect that for detailed dynamics over networks, the lifting operator will create microscopic distributions off away the slow manifold. If so, we can enhance our calculations by forcing the system to start from consistent to the coarse-grained variables (lower-order moments of the microscopic distribution) microscopic initial conditions almost *on* the slow manifold.

Considering the spread of an epidemic in networks, in analogy to spatial models in space [59, 60], the higher order moments are associated with the "spatial" densities of pairs, triples etc. [27].

In particular, let us denote our network by $G(V, E)$, where $V = \{v_i\}, i = 1, 2, ..., N$ is the set of vertices -corresponding to the $N$ individuals-, and $E$ is the set of edges, the links between the individuals. An edge $e_{v_i v_j}$ is defined by $\{v_i, v_j\}$ where $v_i, v_j \in V$ are the nodes associated with it; $e_{v_i v_j} = 1$ if $\{v_i, v_j\}$ are connected and $e_{v_i v_j} = 0$ otherwise. Here all links are bidirectional and self-contacts are not allowed, i.e. $e_{v_i v_i} = 0$.

If $U\{\mathbf{x}_{v_l} = x_i\}$ is the distribution of the discrete epidemic state $x_i = \{S, I, R, ...\}$ over the $L$ vertices $v_l$, $(l = 1, 2, ...., L)$ of the network then its *first* moment denoting the mean value of the state $x_i$ is given by

$$E[x_i] = \frac{1}{L} \sum_{l=1}^{L} \delta_{v_l}(x_i - x_l) \tag{15}$$

where

$$\delta_{v_l}(x_i - x_l) = \begin{cases} 1, & \text{if } U\{\mathbf{x}_{v_l} = x_i\} \\ 0, & \text{otherwise} \end{cases} \tag{16}$$

Since we are dealing with multiple states the second moments are associated with the covariances, reflecting the strength of correlation between two of the epidemic states over the links of the nodes $v_l$. By definition, $Var[x_i] = Cov[x_i, x_i]$ reading:

$$Var[x_i] = \frac{1}{L_{pairs}} \sum_{l=1}^{L} \sum_{e_{v_l v_j}=1} \delta_{v_l}(x_i - x_l) \delta(x_i - x_j) \tag{17}$$

where the sum $\sum_{e_{v_l v_j}=1}$ denotes summing over the nodes $v_j$ that are connected with the node $v_l$; $L_{pairs}$ is defined by $L_{pairs} = \sum_{l=1}^{L} \sum_{e_{v_l v_j}} 1$ and corresponds to the number of pairs in the network. If the connectivity degree is a constant number for all nodes, say $c$, then the above sum gives $L_{pairs} = c L$. Note that according to (15), (17) the pairs are counted twice. Hence the correlation sum will always result to an even number (see figure 3a for a simple example).

The other elements of the covariance matrix are given by

$$Cov[x_i, x_j] = \frac{1}{L_{pairs}} \sum_{l=1}^{L} \sum_{e_{v_l v_k}=1} \delta_{v_l}(x_i - x_l) \delta(x_j - x_k) \tag{18}$$

In the case of a network, the definition of the third moment is not unique: triplets can arise as chain-like (figure 3a) or loop-like connections (figure 3b).

**2.2.2 The computational procedure: Coupling the Equation-Free approach with Simulated Annealing**

The convergence to the slow manifold at specific initial values of the coarse-grained slow variables can be achieved by exploiting the Equation-Free approach as follows (see figure 4):

(1) Prescribe the desired coarse-grained initial conditions $x(t_0)$.

(2) Transform $x(t_0)$ through a lifting operator $\mu$ to consistent microscopic realizations: $U(t_0) = \mu\, x(t_0)$. In general, the higher order moments of the microscopic distributions, say $y[U(t_0)]$, will be off the slow manifold.

*k=0;*

*Do until convergence to the slow manifold {*
  *k=k+1;*

(3) Evolve these realizations in time using the microscopic simulator for a very short macroscopic time $dT \ll T$, generating the value(s) $U^k(t_0 + dT)$ and compute the higher moments (for practical means, up to a specific order *l*) $y[U^k(t_0 + dT)]$.

(4) Restrict back to the prescribed coarse-grained initial conditions $u_s(t_0)$, preserving the values of $y[U^k(t_0 + dT)]$.

At this point, Simulated Annealing [40-42] is exploited in order to design the micro-structure of the network.

The objective function at the step $j$ of the SA algorithm may be defined as

$$O^j = \|y^j - y(U^k)\|, \tag{19}$$

where $y(U^k)$ are the target values of the fast variables. In particular within a network with discrete-value coarse-grained states (e.g. susceptible, infected, recovered, isolated etc), the fast variables $y$ correspond to the densities of pairs, triples, quadruples etc. links of individuals of certain states in the network which are related to the second, third, and fourth spatial moments of the underlying distribution respectively [18, 60].

In this case, the proposed SA algorithm reads as follows:

*0) Set the initial systems pseudo-temperature Temp and select the annealing schedule, i.e. the way the pseudo-temperature will decrease.*

*Do until convergence {*
  *1) Evaluate the pseudo-energy (objective function) $O(y)$ of the network;*

  *2a) Select randomly an individual (or a set of individuals)*
  *2b) Select randomly another individual (or a set of individuals)*
  *2c) Create a new state distribution over the network by interchanging the states of the -selected from (2a) and (2b)- individuals(s). The interchange (e.g. between a susceptible and an infected and vice versa) leaves the first-order moment unaltered (see figure 5).*

  *Let us denote the new values of the fast variables corresponding to the new network configuration by $y'$. Then:*

  *3) Evaluate the new objective function $O(y')$.*

*4) Accept or reject the new network configuration using the Metropolis procedure* [61]:

*4a) Accept the new configuration if $O(y') < O(y')$,*

*4b) Accept the new configuration if $O(y') > O(y')$, with a probability*

$$exp\left[-(O(y) - O(y'))/Temp\right];$$

*Otherwise, reject it.*

*5) Reduce the system pseudo-temperature according to the annealing schedule.*

*} End Do*

*} End Do*

### 3. The case study: a simple Individual-Based Epidemic Model on a RRN

Here, the epidemic evolves in discrete time $t$ on a regular random graph, which is characterized by a constant connectivity $d$ between individuals, here, set equal to four. The graph simulates the social structure of our artificial individual-oriented world involving $N$ individuals. Here in contrast with lattice-based models such as Cellular Automata the spatial position of the individuals is rather irrelevant: individuals are connected in a random, non-local way which resembles better the social interaction mechanism of such problems [62]. The links are bidirectional and fixed during the epidemic and loops (connections of individuals to themselves) are not allowed.

Each individual is labeled by an index $i = 1, 2, ..., N$ and is characterized just by its health state respect to the epidemic: (a) Susceptible ($S$) when the individual is not yet infected but there is a certain probabilistic potential to get infected; (b) Infected ($I$) when is a carrier of the disease and can potentially transmit it to its links and (c) Recovered ($R$) when the individual recovers from the infection, cannot transmit and is temporally immunized from the infection. For the calculations the state of each individual is coded with the indices 1, 2 or 3 according to whether the individual is a susceptible, infected or a recovered one.

The random regular graph was generated based on the procedure described in [63]:

1. An unlinked pair $(i, j)$ of individuals with $i \neq j$ is randomly chosen from a total of $Nd$ even number of points ($d$ points in $N$ groups).

2. Connect $i$ with $j$ setting $G_{ij} = 1$; leave $G_{ij} = 0$ otherwise, and

3. Repeat the above procedure until all individuals are properly linked and the graph is complete, i.e. the infection can potentially reach every individual in the network when starting from any other one.

Other approaches for generating asymptotically random regular graphs for larger degrees of connectivity can be found in [62].

In our artificial world, individuals interact in the network with their links and change their states day by day in a probabilistic manner according to simple rules involving their own states and the states of their links.

At each discrete time step $t$ of the simulation, the following rules are updated in a synchronous way:

- Rule #1: An infected individual $i$ infects a susceptible link $j$ with a probability $p_{S \to I}$.

- Rule #2: An infected $i$ individual recovers with a probability

$$p_{I \to R} = 1 - \exp\left(-c_1 f(i)^{-c_2}\right) \qquad (23a)$$

$f(i)$ is the ratio of infected links of the infected individual:

$$f(i) = \frac{n_i(I)}{d} \qquad (23b)$$

where $n_i(I)$ denotes the number of infected links of the infected individual $i$ at time $t$.

The above relation reflects the fact that the probability of passing into the recovery phase depends in a nonlinear way on the number of infected neighbors: the recovery from infection becomes more difficult for infected individuals surrounded by many infected neighbors and thus the infection persists for a longer time. The nonlinearity may be accounted to different factors such as a drift of the disease-virus over short time periods and heterogeneity in infectiousness and/or recovery. Actually, the motivation for using relation (1) was an analogous empirical one proposed in [12] and modified in [43] giving the average probability that an individual fails to clear an infection developing the carrier state. Over the past years have been considered various other forms of nonlinear infection incidence and recovery rates [64-66]. The parameters $c_1$ and $c_2$ may be interpreted as parameters dependent from the disease characteristics and may be determined by statistical data.

- Rule #3: A recovered individual becomes susceptible again with probability $p_{R \to S}$ otherwise cannot transmit and possesses temporal immunity against the disease. Here, the probability $p_{R \to S}$ represents the reciprocal expected recovery period of the disease.

## 4. Simulation Results and Discussion: Coarse-Grained Bifurcation and Rare-Event Computations

The purpose, here, is to perform systematic coarse-grained tasks such as bifurcation and stability analysis of the epidemic model without extracting analytically any closed macroscopic equations such as mean-field or pair-wise approximations. Instead, the individual-based epidemic model is exploited using the Equation-free approach as an input-output coarse "black-box" timestepper.

Simulation results were obtained on a fixed RRN of $N =20,000$ individuals each and every one having $d = 4$ links while the values of the parameters determining the transition probability from the infective to the recovered state $p_{I \to R}$ and the transition probability from the recovered to the susceptible state $p_{R \to S}$, were set as $c_1 = 0.1$ and $c_2 = 0.5$ and $p_{R \to S} = 0.2$.

Figure 6 shows some characteristic temporal simulations for a wide range of values of the parameter $p_{S \to I}$. It is clearly shown that depending on the value of $p_{S \to I}$, the dynamics of the epidemic model reveal some interesting nonlinear behavior. In particular, for relatively big values of the $p_{S \to I}$ there is only one stable stationary solution corresponding to a high endemic situation (figure 6a, b). The same behaviour is observed for even smaller values of $p_{S \to I}$ (figure 6c).

When the value of $p_{S \to I}$ decreases further and now depending on the initial conditions, a second stable stationary point of low endemicity appears (figure 6d). With a further decrease of the infection transmission probability, it seems that the state of high endemicity disappears giving its place to a low-endemicity situation.

For even smaller values of $p_{S \to I}$ the network converges to the "disease-free" state (figure 6f). These observations suggest that the system has some critical points in the parameter space which mark the onset of these phenomena and -as the bifurcation theory dictates - the existence of unstable coarse-grained states which are unreachable through plain long-run temporal simulations.

Based on the above simulation results, there are two basic issues that should be resolved here, namely, (a) construct the complete set of fixed point solution branches (both stable and unstable ones), (b) locate the values of the parameters for which these critical phenomena emerge.

Both questions are being answered by performing coarse-grained bifurcation analysis exploiting the proposed procedure.

As described in the previous section, the first step in the methodology refers to the selection of the appropriate coarse-grained statistics and the lifting operator. The assumption is that the higher-order moments evolve quickly to a low-dimensional slow manifold parametrized by the lower-order ones. That is, after a relatively "short" macroscopic time-period, higher-order moments are becoming functions of the lower-order moments.

Along the lines of the mean-field approximation one could choose the densities of the susceptible, infected or recovered as the potential coarse-grained (observable) variables and as a lifting operator, the intuitively for the structure of our network simplest one: random susceptible and infection state assignments over the network consistent with the corresponding densities. The hypothesis here is that coarse-grained dynamic behavior of the epidemic spread lies on a low-dimensional manifold parametrized only by these two coarse-grained observables. But in what level do we have a good representation of the coarse-grained dynamics?

Figure 7 depicts the 3-D phase portraits $[S], [I]$ vs. the $[SI]$ pair density around the stationary state at $p_{S \to I}$ = 0.17. Microscopic state configurations consistent with the first and second-order moments over the network were created by utilizing the proposed SA procedure initializing these coarse-grained quantities at will. As it is clearly shown, there are two time-scales at this representation: (i) an initial fast approach to an one-dimensional manifold which can be parametrized by either $[S]$ or $[I]$; and (ii) a slow approach to the ultimate stationary state - at $([S], [I]) = (0.034, 0.903)$- on this slow manifold.
This implies that the expected dynamics, over the ensembles conditioned on this slow manifold, can "close" in terms of the first-order moments.

As explained in the previous section, the atomistic simulator was treated as a black-box coarse-grained timestepper of the macroscopic –observed variables. The time-horizon was set as $T = 4$ time steps (days) while the parameter $dT$ in the SA procedure was set as $dT = 1$ time step.

Figure 8 shows the resulting one-parameter coarse-grained bifurcation diagram of $[I]$ with respect to $p_{S \to I}$. There are two regular turning points at $p_{S \to I} \approx 0.138$ and $p_{S \to I} \approx 0.15$ that mark the barrier between coarse-grained stable (solid lines) and unstable (dotted lines) equilibria. Continuation around the turning points is accomplished by coupling the fixed point algorithm with the pseudo-arclength continuation technique. The coarse-grained solutions found by the proposed approach are consistent with long-run temporal simulations. The diagram of the algebraically largest eigenvalue of $\frac{\partial \boldsymbol{\Phi}_T}{\partial \boldsymbol{x}}$ with respect to the bifurcation parameter is also shown in the inset of figure 8.

## 5. Conclusions

One of the most important issues in mathematical epidemiology revolves around the investigation of the emergent dynamics at the macroscopic scale. Traditionally, for the systematic study of the dynamics one has to resort to numerical bifurcation analysis tools in order to trace solution branches through bifurcation points, locate bifurcation sets and critical points in parameter space, identify the regions of hysteresis phenomena, etc. A fundamental prerequisite for such tasks is the availability of reasonably accurate closed form dynamical models. However real-world epidemiological problems are characterized –due to their stochastic/microscopic nature and nonlinear complexity- by the lack of such good explicit, coarse-grained macroscopic evolution equations.

Detailed individual (agent-based)-based multiscale models are the state-of-the-art in the field and are often used as they are considered to incorporate the appropriate level of complexity. When this is the case, conventional continuum algorithms cannot be used directly. Severe problems arise in trying to find the closures to bridge the gap between the scale of the available description and the macroscopic scale at which the questions of interest are asked and the answers are required. Hence, what is usually done with such complex, detailed epidemic simulators is simple simulation: initial conditions and operating conditions are set, and then the simulation is "run" (evolved) over time. To understand the dynamic behaviour, many simulations have to be performed by changing initial conditions, parameters, and interaction rules. Yet, simulation forward in time is but the first of systems level computational tasks one wants to explore while analyzing the parametric behavior of an epidemiological model; it is both time consuming and for several tasks inappropriate (e.g. unstable solutions cannot be found).

In this work it is shown how the Equation-Free methodology for multiscale computations can be exploited to systematically analyse the coarse-grained dynamics of individualistic epidemic simulators on networks under the pairwise correlations perspective. A Simulated Annealing procedure is proposed in order to drive the coarse timestepper to its own coarse-grained slow manifold on which the coarse-grained dynamics evolve after a fast-in the macroscopic scale- time interval. Using this framework steady state and stability computations, continuation and numerical bifurcation analysis of the complex-emergent dynamics can be performed in a full-computational manner bypassing the need of analytical derivation of closures for the macroscopic-level equations.

For illustration purposes I used an extremely simple individual-based epidemic model deploying on a caricature of a social network. While the simulation example is admittedly simple, it does demonstrate the scope of the system level tasks that one can attempt using the proposed framework by acting directly on the individual-based epidemic simulation on networks. Stability and bifurcation analysis have been demonstrated here; more complex and closer to the real-world social networks, rare-events, coarse projective integration for the acceleration in time of the simulation as well as coarse-grained model-predictive, nonlinear control design and optimization can also be attempted.

Efforts could proceed along a wide directions ranging from the analysis of the dynamic interplay between network state and network topology to the systematic reconstruction of the slow manifold using advanced techniques such as the Computational Singular Perturbation method for the [67,68]. The use of state-of-the-art data mining techniques such as Diffusion Maps [69, 70] that can be exploited to efficiently extract the "correct" coarse-grained variables from a more complex individual-based large-scale code could also be attempted.

**Figures**

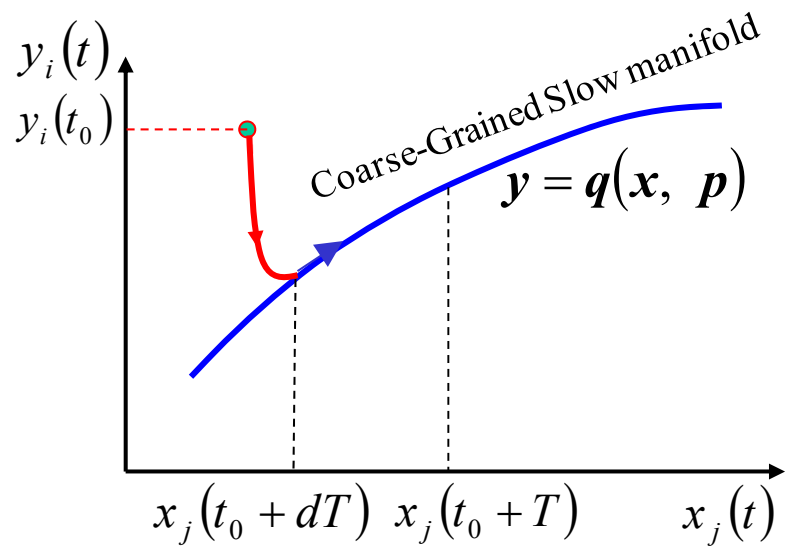

**Figure 1.** The basic assumption of the Equation-Free approach: Very fast the dynamics of the complex systems evolve on a slow coarse-grained manifold.

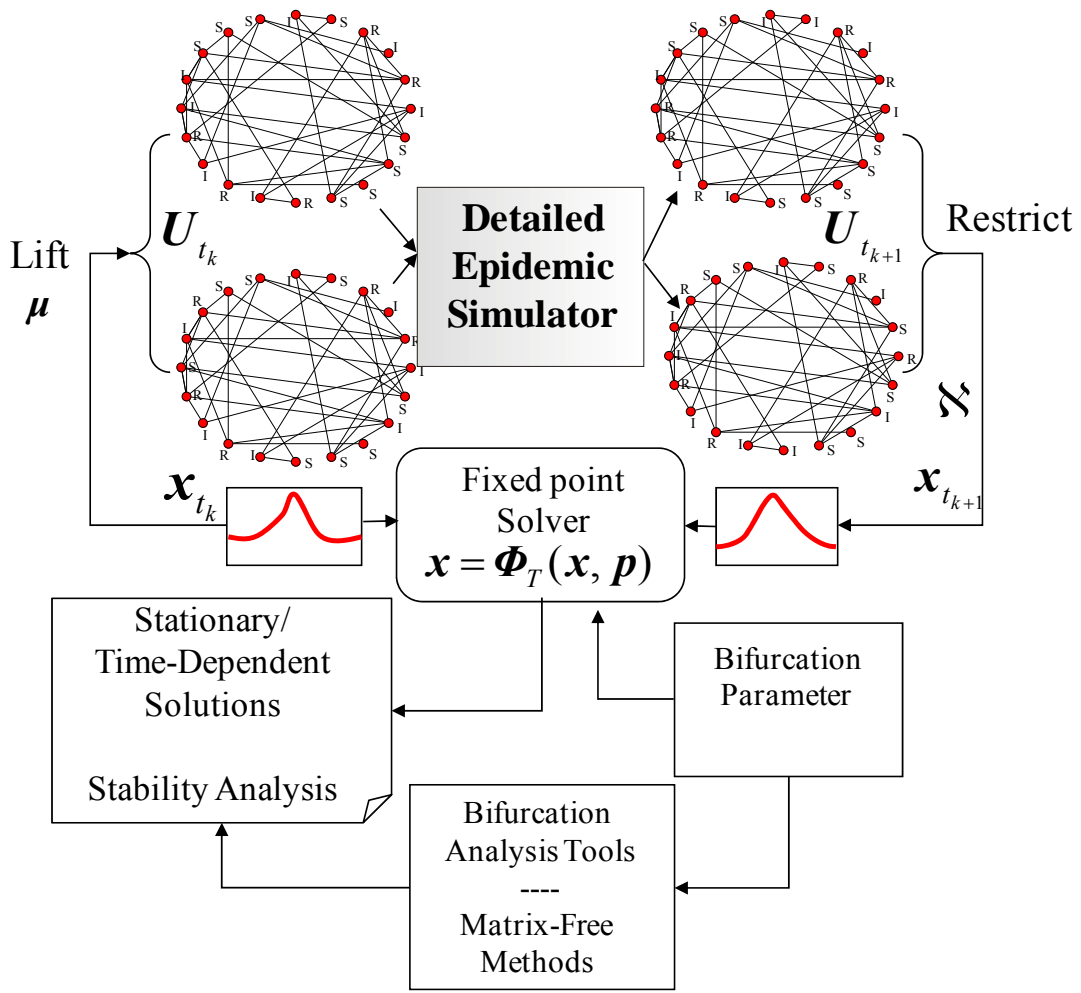

**Figure 2.** Schematic of the Coarse Timestepper Concept.

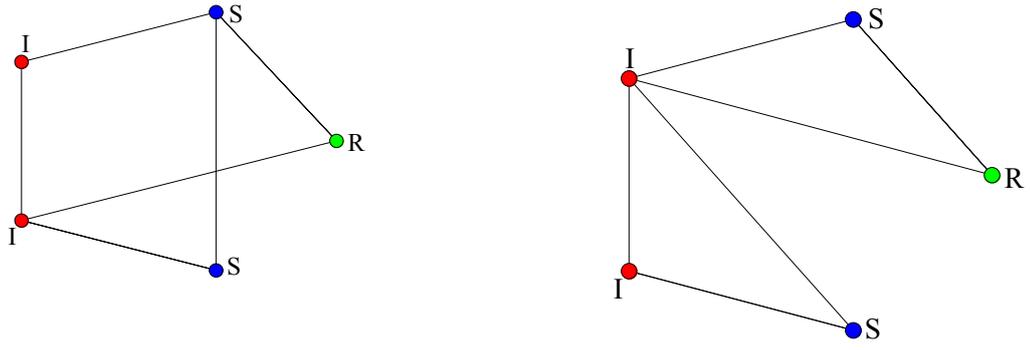

**Figure 3**. Simple examples with a network of 5 nodes with (a) 12 pairs: [IS], [SI], [SR], [RS], [SS], [SS], [SI], [IS], [IR], [RI], [II], [II] and 12 chain-like triples: [IIS], [SII], [SIR], [RIS], [ISS], [SSI], [IRS], [SRI], [SSR], [RSS], [ISR], [RSI], (b) 12 pairs [IS], [SI], [SR], [RS], [RI], [IR], [IS], [SI], [II], [II], [IS], [SI], 12 loop-like triples: [ISR], [SRI], [RIS], [RSI], [SIR], [IRS], [IIS], [ISI], [SII], [IIS], [ISI], [SII], and 10 chain-like triples: [IIS], [SII], [SIR], [RIS], [SIS], [SIS], [IIR], [RII], [RII], [IIR].

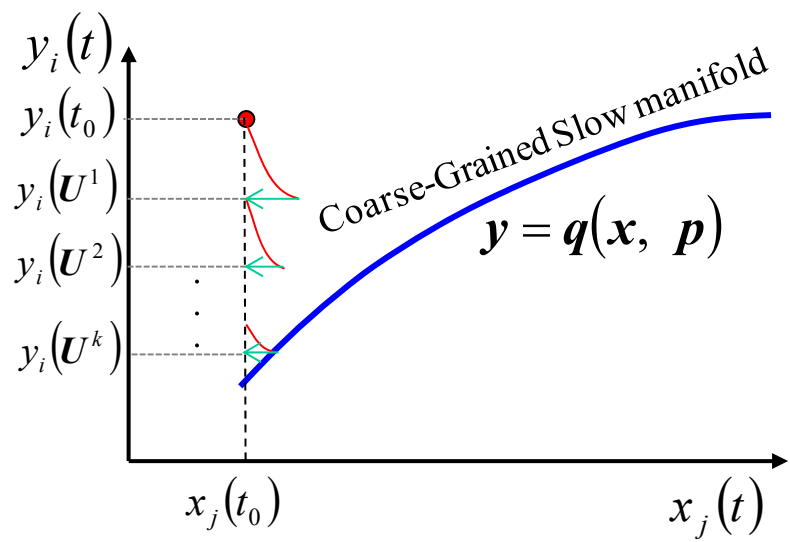

**Figure 4**. Using Simulated Annealing to drive the system on the attracting coarse-grained slow manifold.

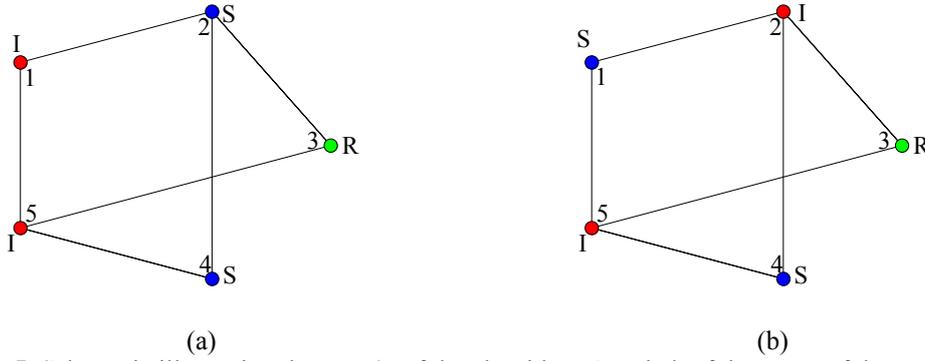

(a)                          (b)

**Figure 5**. Schematic illustrating the *step 2c* of the algorithm. A switch of the states of the nodes 1&2 does not affect the mean densities: $[S]=\frac{2}{5}, [I]=\frac{2}{5}, [R]=\frac{1}{5}$ but changes the microstructure, i.e. the higher order moments of the network; for example the density of the pairs change from

(a) $[SS]=\frac{2}{12}, [SI]=\frac{4}{12}, [SR]=\frac{2}{12}, [II]=\frac{2}{12}, [IR]=\frac{2}{12}$ to (b) $[SI]=\frac{8}{12}, [IR]=\frac{4}{12}$

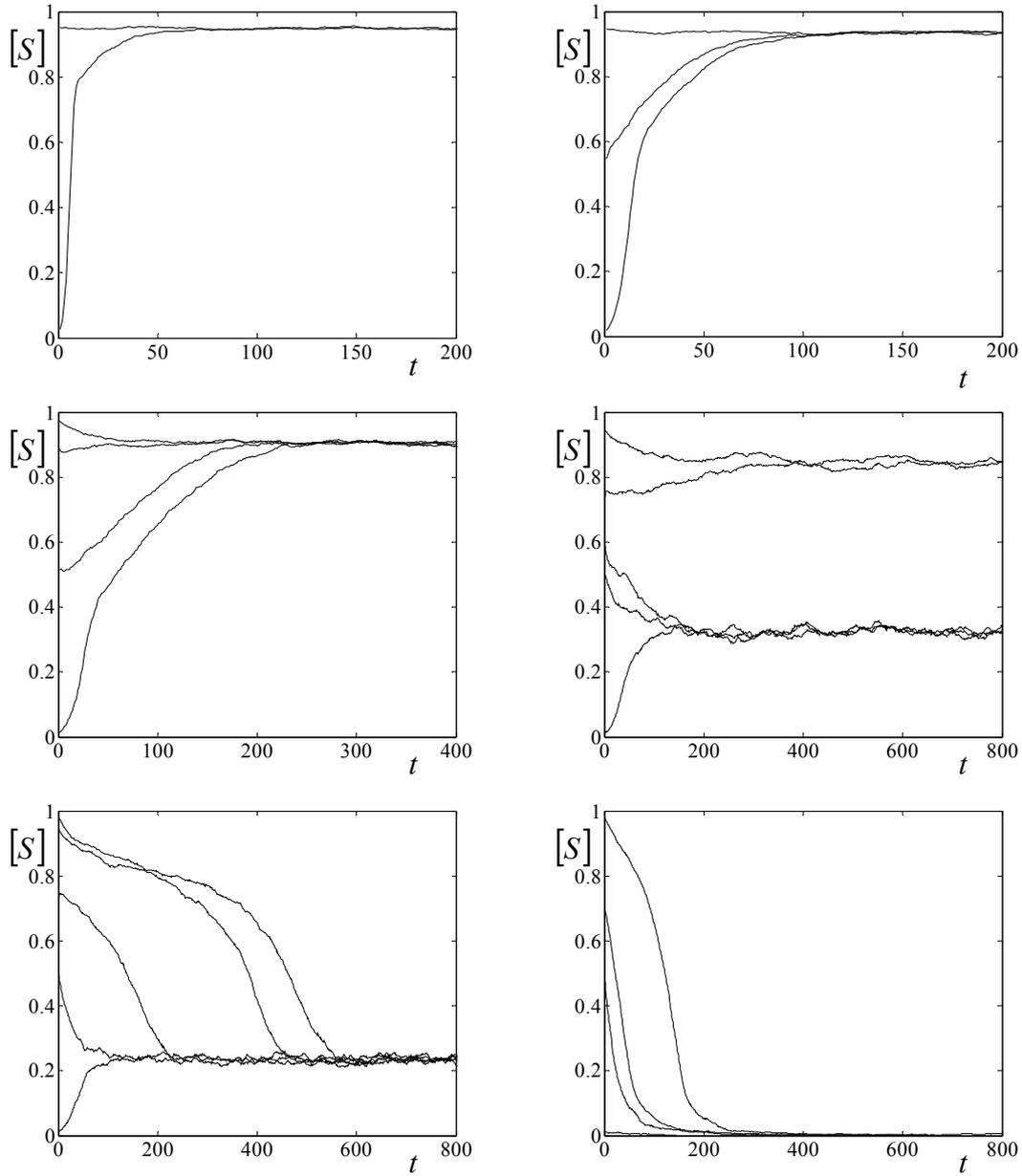

**Figure 6.** Temporal simulation showing the evolution of $[S]$ of the individual-based epidemic simulator with $c_1 = 0.1$ and $c_2 = 0.5$ and $p_{R \to S} = 0.2$. (a) $p_{S \to I} = 0.5$, (b) $p_{S \to I} = 0.25$, (c) $p_{S \to I} = 0.17$, (d) $p_{S \to I} = 0.14$, (e) $p_{S \to I} = 0.13$, (f) $p_{S \to I} = 0.10$

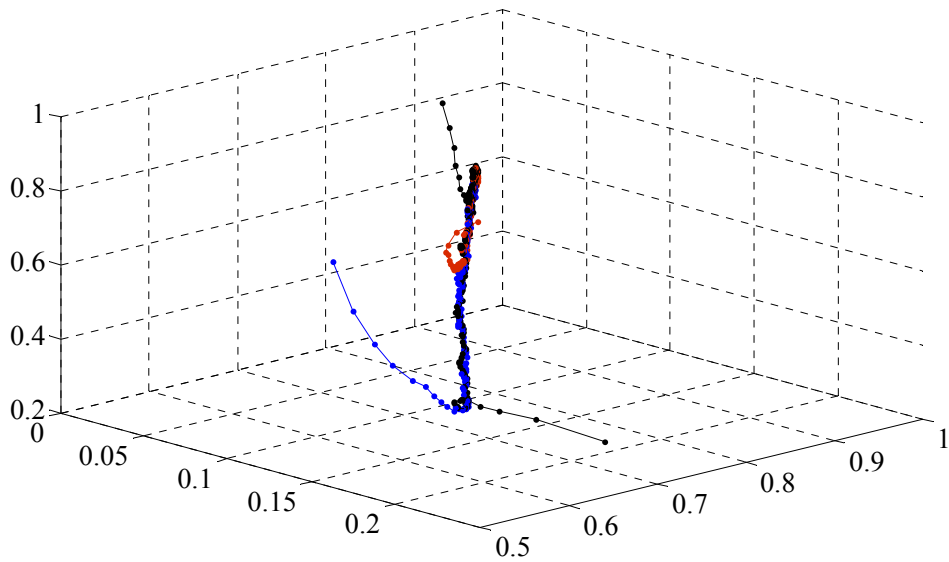

**Figure 7.** Phase portrait of $[S],[I]$ vs. the $[SI]$ pair density around the stationary point at $p_{s \to I}$ = 0.17. Initial conditions of the pair densities were constructed at will using the proposed SA scheme. An one-dimensional manifold is clearly shown.

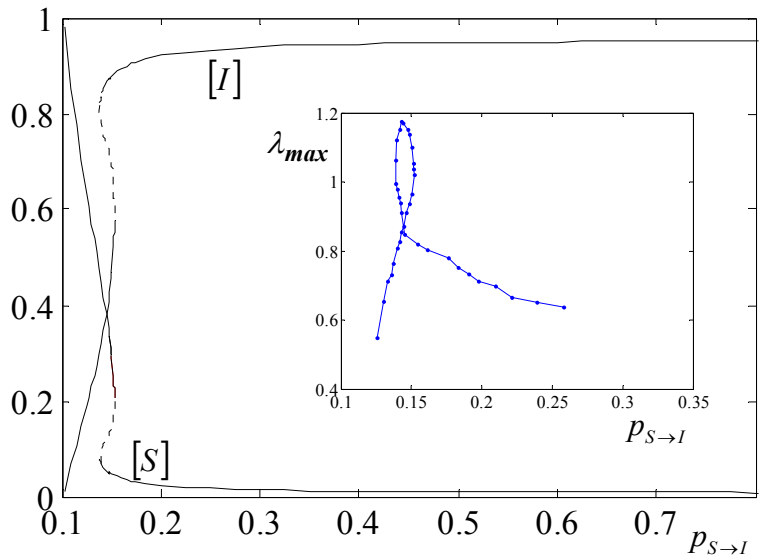

**Figure 8.** Coarse-grained bifurcation diagrams of $[I]$ and $[S]$ vs. $p_{S \rightarrow I}$ as obtained using the microscopic simulator as a black box coarse timestepper. Solid lines depict coarse-grained stable states while dotted lines depict coarse-grained unstable states. The inset shows the diagram of the algebraically largest eigenvalue with respect to the bifurcation parameter.